\documentclass{siamltex}		
\usepackage[utf8]{inputenc}
\usepackage{graphicx,multirow}  
\usepackage{verbatim,float} 
\usepackage{todonotes}
\usepackage{amssymb,amsmath}
\usepackage[T1]{fontenc}
\usepackage{graphicx,amsmath,amsfonts,amssymb,subfigure}
\usepackage{color}
\usepackage{hyperref}
\usepackage[mathscr]{euscript}
\usepackage{algorithm,algorithmic}

\newcommand{\eq}[1]{\begin{equation}\label{#1}}
\newcommand{\en}{\end{equation}}
\usepackage{multirow}

\graphicspath{{./FIGS/}}

\usepackage{verbatim}

\title{RBF approximation of  three dimensional PDEs using  Tensor Krylov subspace methods }
\author{M. El Guide \thanks{Centre for Behavioral Economics and Decision Making(CBED), FGSES, Mohammed VI Polytechnic University, Green City, Morocco}   \and K. Jbilou\footnotemark[3] \thanks{LMPA, 50 rue F. Buisson, ULCO Calais, France; Mohammed VI Polytechnic University, Green City, Morocco; jbilou@univ-littoral.fr }  \and A. Ratnani\thanks{Mohammed VI Polytechnic University, Green City, Morocco}}

\begin{document} 
	
	\maketitle 
	\begin{abstract}
		In this paper, we propose different algorithms for the solution of a tensor linear discrete ill-posed problem arising in the application of the meshless  method for solving  PDEs in three-dimensional space using multiquadric radial basis functions. It is well known that the truncated singular value decomposition (TSVD) is the most common effective solver for ill-conditioned systems, but unfortunately the operation count for solving a linear system with the TSVD is computationally expensive for large-scale matrices. In the present  work, we  propose algorithms based on the use of the well known Einstein product for two tensors to define the tensor global Arnoldi and the tensor Gloub Kahan bidiagonalization algorithms.  Using the so-called Tikhonov regularization technique, we will be able to provide computable approximate
		regularized solutions in a few iterations.
	\end{abstract}
	
	\begin{keywords} 
		Krylov subspaces, Linear tensor equations,  Tensors, Global Arnoldi, Global Golub-Kahan, Einstein product.
	\end{keywords}
	

	\section{Background and introduction}
	The most commonly used numerical methods for solving partial differential equations are finite element, finite difference and finite volume. However, a lot of work has to be done to generate meshes when using these methods and the task becomes more difficult  when  dealing  with
	complicated domains in  higher dimensions, i.e., $d \geq  3$. That is when  meshless methods \cite{Liu} comes into the picture. The  idea of
	applying  the meshless methods  in many  branches  of science  and engineering  has
	gained popularity over  the past few years,  e.g., in elasticity,  ground and water  flow, wave
	propagation  and in  electromagnetic  problems \cite{GW, WTBH}. Meshless methods are able to deal with any type of PDE and in any space dimension, by only using of the so-called radial basis functions (RBF). With their radial symmetries properties, RBF  demonstrate the ability to transform a problem in several dimensions into one that is one-dimensional. The greatest advantage of these methods is the possibility of enriching the functions, that is to say,  including physics properties of the studied problem. We can
	thus, with a small number of nodes, approach with great precision the solution of the problem. Among these meshless techniques, the Multiquadric (MQ) proposed by Rolland Hardy in \cite{Hardy}  is the most commonly RBF used in applications. Other common used radial basis function
	can be consulted in \cite{Fass, Wend}. The wide use of MQ-RBF  is due to its global,
	infinitely differentiable properties which make it a good candidate to give good approximation properties. In this paper we will focus on the use of MQ as RBF approximation for solving three dimensional PDEs, which will result in a multidimensional linear equation. This comes at a cost since the matrix problem
	to solve in the approximation becomes dense. The use of MQ will give rise to an ill-posed linear system of equations, i.e. the system matrix is critically
	conditioned, which leads to numerical implementation difficulties. Although a finite element discretization of the
	problem yields a sparse well conditioned matrices, it
	requires a discretization of the whole domain 
	which is not
	always feasible, e.g., if the domain is unbounded. Though
	the topic of solving ill-conditioned and dense linear systems built upon meshless methods has been around for a number of years, the lack of
	efficient linear solvers and the little attention that has been focused  on the conditioning of the linear systems has delayed a full exploration of meshless–based approaches.  Regularization techniques are then needed to reduce the effect of the ill-conditioning of the system matrix. When using a large number  of MQ-RBF approximation points, the problem of the number of
	floating point operations needed to solve a linear system should be also addressed. It is the purpose of this paper to overcome the aforementioned
	difficulties by combining regularization techniques and iterative methods based on Krylov subspace techniques. It's worth to mention that Krylov methods such as GMRES and LSQR have been already used  for meshless methods using radial basis functions; see, e.g., \cite{BCM}. In this paper a meshless method based on multiquadric (MQ) radial basis function is proposed for solving numerically the modified Helmholtz equation.
	
	This paper is organized as follows.  We shall first   present in  Section 2 
	some symbols and notations used throughout the paper. We also recall the concept of contract product between two tensors. In Section 3 we describe the multiquadric radial basis function framework and then we propose  a tensor formulation of the discretization of the three-dimensional PDEs. In Section
	4, we present some inexpensive approaches based on global Krylov subspace methods combined with regularization techniques to solve the obtained
	ill-posed dense linear systems arising when using multiquadric radial basis function. Section 5 is
	dedicated to some numerical experiments. Concluding remarks can be found in Section 6.
	
	\section{Definitions and Notations}\label{sec:section2}   In this section, we briefly review some concepts and notions
	that are used throughout the paper. 	A tensor is a multidimensional array of data and a natural extension of scalars, vectors and matrices to a higher order. Notice that a  scalar is a $0^{th}$ order tensor, a vector is a  $1^{th}$ order tensor and a matrix is   $2^{th}$  order tensor. The tensor order  is  the number of its indices, which is called modes or ways. For a given N-mode tensor $  \mathcal {X}\in \mathbb{R}^{I_{1}\times I_{2}\times I_{3}\ldots \times I_{N}}$,
	the notation $x_{i_{1},\ldots,i_{N}}$ (with $1\leq i_{j}\leq I_{j},\; j=1,\ldots N $) stand for element $\left(i_{1},\ldots,i_{N} \right) $ of the tensor $\mathcal {X}$. Corresponding to a given tensor $ \mathcal {X}\in \mathbb{R}^{I_{1}\times I_{2}\times I_{3}\ldots \times I_{N}}$, the notation $$  \mathcal {X}_{\underbrace{:,:,\ldots,:}_{(N-1)-times},k} \quad k=1,2,\ldots,I_{N}$$ denotes a tensor in $\mathbb{R}^{I_{1}\times I_{2}\times I_{3}\ldots \times I_{N-1}}$ which is obtained by fixing the last index and is called frontal
	slice. Fibers are the higher-order analogue of matrix rows and columns. A fiber is
	defined by fixing every index but one. A matrix column is a mode-1 fiber and a
	matrix row is a mode-2 fiber. Third-order tensors have column, row, and tube fibers; see  \cite{kimler1,kolda1} for more detail . Throughout this work, vectors and matrices are respectively denoted by
	lowercase and capital letters, and tensors of  higher order  are represented by  
	calligraphic   letters.
	
	\noindent We first recall the definition of the well known  $n$-mode tensor product; see \cite{kolda1} .
	
	\begin{definition}
		The $n$-mode product of the tensor 	$\mathcal {A}=[a_{i_1i_2 \ldots i_n} ]  \in \mathbb{R}^{I_{1} \times I_2\times \ldots \times I_N}$ and the matrix $U=[u_{j,i_n}] \in \mathbb{R}^{J \times I_n}$ denoted by $\mathcal {A} \times_n U$ is a tensor of order $I_1\times I_2 \times \ldots \times I_{n-1} \times J \times I_{n+1} \times \ldots \times I_N$ and defined by
		\begin{equation}
		\label{nmode}
		(\mathcal {A} \times_n U)_{i_1i_2\ldots i_{n-1}ji_{n+1}\ldots i_N}= \displaystyle \sum_{i_n=1}^{I_N} a_{i_1i_2\ldots i_N} u_{j,i_n}
		\end{equation}
	\end{definition}
	
	\medskip 
	
	\noindent The $n$-mode product of the tensor $\mathcal {A} \in \mathbb{R}^{I_{1} \times I_2\times \ldots \times I_N}$ with the vector $v=[v_{i_n}] \in \mathbb{R}^{I_n}$  is denoted by 
	$\mathcal {A} \bar {\times} v$ and given by
	$$ (\mathcal {A} \bar {\times} v )_{i_1\ldots i_{n-1}i_{n+1}\ldots i_N}= \displaystyle \sum_{i_n=1}{I_n}x_{i_1i_2\ldots i_N} v_{i_n} .$$ 
	
	\medskip 
	
	\noindent For this $n$-product, we have the following properties. Let $\mathcal {A}= \in \mathbb{R}^{I_{1} \times I_2\times \ldots \times I_N}$ and consider two matrices $B \in  \mathbb{R}^{J_m \times I_m}$ and $C  \in  \mathbb{R}^{J_ n \times I_n}$ with $m \ne n$. Then
	$$A \times_n B \times_m C=A \times_m C \times_n B$$
	and if $B \in  \mathbb{R}^{J \times I_n}$ and $C  \in  \mathbb{R}^{I_ n \times J}$ , then
	$$ A \times_n B \times_m C=A \times_n (CB).$$

	\noindent 	Next, we recall the definition and some properties  of the tensor Einstein product which is an  extension of the matrix product; for more details see \cite{brazell}
	
	\begin{definition}\cite{einstein}
		
		Let $\mathcal {A} \in \mathbb{R}^{I_{1}\times I_{2}\times \ldots \times I_{L}\times   K_{1}\times K_{2}\times \ldots \times K_{N}} $, $\mathcal {B}, \in \mathbb{R}^{K_{1}\times K_{2}\times \ldots \times K_{N}\times J_{1}\times J_{2}\times \ldots \times J_{M}} $, the Einstein product of tensors $\mathcal {A}$ and $\mathcal {B}$ is a tensor of size   $  \mathbb{R}^{I_{1}\times I_{2}\times \ldots \times I_{L}\times J_{1}\times J_{2}\times \ldots \times J_{M}} $ defined as  :
		$$ (\mathcal {A}\ast_{N}\mathcal {B})_{i_{1}\ldots i_{L} j_{1} \ldots j_{M}}=       \sum_{k_{1}=1}^{K_{1}}\sum_{k_{2}=1}^
		{K_{2}}\sum_{k_{3}=1}^{K_{3}}\ldots \sum_{k_{N}=1}^{K_{N}}a_{i_{1}\ldots i_{L} k_{1} \ldots k_{N} }b_{k_{1} \ldots k_{N}j_{1} \ldots j_{M}}.  $$
		\label{einstein}
	\end{definition} 
	
	\noindent Let $\mathcal {A} \in \mathbb{R}^{I_{1}\times I_{2}\times \ldots \times I_{N}\times   J_{1}\times J_{2}\times \ldots \times J_{M}} $ and let $\mathcal {B} \in \mathbb{R}^{J_{1}\times J_{2}\times \ldots \times J_{M}\times   I_{1}\times I_{2}\times \ldots \times I_{N}} $ such that $b_{i_1\ldots i_Mj_1 \ldots j_m}= a_{j_1 \ldots j_Ni_1 \ldots i_M}$. Then $\mathcal {B} $ is called the transpose of  $\mathcal {A}$ and denoted by $\mathcal {A}^T$.

	\noindent  A tensor $\mathcal{D}=[d_{i_{1},\ldots,i_{M},j_{1},\ldots,j_{N}}] \in \mathbb{R}^{I_{1}\times I_{2}\times \ldots \times I_{N}\times   I_{1}\times I_{2}\times \ldots \times I_{N}} $ is a diagonal tensor if $d_{i_{1},\ldots,i_{N},j_{1},\ldots,j_{N}}=0$ in the case that the indices $ i_{1},\ldots,i_{N} $ are different from $j_{1},\ldots,j_{N}$. If $\mathcal{D}$ is a diagonal tensor such that all the diagonal entries are equal to $1$, then $\mathcal{D}$ is the unit tensor denoted by $\mathcal{I}_N$. 
	$\mathcal{O}$ is the tensor having all its entries equal to zero.. \\
	
	\begin{definition}
		Let  $\mathcal{A} \in \mathbb{R}^{I_{1}\times I_{2}\times \ldots \times I_{N}\times   I_{1}\times I_{2}\times \ldots \times I_{N}} $. The tensor  $\mathcal{A}$ is invertible if there exists a tensor $\mathcal {X} \in  \mathbb{R}^{I_{1}\times I_{2}\times \ldots \times I_{N}\times   I_{1}\times I_{2}\times \ldots \times I_{N}} $ such that
		$$\mathcal{A}\ast_{N}\mathcal {X}=\mathcal{X}\ast_{N}\mathcal {A}=\mathcal {I}_N$$
		where $\mathcal {I}_N$ denotes the identity tensor. In that case, $\mathcal{X}$ is the inverse of $\mathcal{A}$ and is denoted by  $\mathcal{A}^{-1}$.\\
	\end{definition}
	
	\medskip
	\noindent The trace of an even-order  tensor $\mathcal{A}\in \mathbb{R}^{I_{1}\times I_{2}\times I_{3}\ldots \times I_{N}\times   I_{1}\times I_{2}\times I_{3}\ldots \times I_{N}}$ is given by 
	\begin{equation}
	\label{tr1}
	tr(\mathcal{A})=\sum_{i_{1} \ldots i_{N}} a_{i_{1} \ldots i_{N} i_{1} \ldots i_{N}}.
	\end{equation}
	\medskip
	\noindent We have the following relation. Let $ \mathcal{A}\in \mathbb{R}^{I_{1}\times I_{2}\times I_{3}\ldots \times I_{N}\times   J_{1}\times J_{2}\times \ldots \times J_{M}}$ and $\mathcal{B}\in \mathbb{R}^{J_{1}\times J_{2}\times \ldots \times J_{M}\times   I_{1}\times I_{2}\times I_{3}\ldots \times I_{N}}$, then
	\begin{equation}
	\label{tr2}
	tr(\mathcal{A} \ast_M \mathcal{B})= tr(\mathcal{B} \ast_N \mathcal{A})
	\end{equation}
	
	\begin{definition}
		The inner product of two tensors   of
		the same size  $\mathcal {X},\mathcal {Y}\in \mathbb{R}^{I_{1}\times \ldots \times I_{N}\times   J_{1}\times \ldots \times J_{M}}$
		is given by: 
		\begin{equation}\label{inner1}
		\langle \mathcal {X},\mathcal {Y} \rangle=tr(\mathcal {X}^{T}\ast_{N}\mathcal {Y}) 
		\end{equation}
		where $\mathcal {Y}^{T}\in \mathbb{R}^{J_{1}\times J_{2}\times J_{3}\ldots \times J_{M}\times I_{1}\times I_{2}\times I_{3}\ldots \times I_{N}}$ denote de transpose of $\mathcal {Y}.$ \\
		The Frobenius norm of the tensor  $mathcal {X}$ is given by 	
		\begin{equation}\label{inorm11}
		||\mathcal {X}||_F=\displaystyle \sqrt{tr(\mathcal {X}^{T}\ast_{N}\mathcal {X}) }.
		\end{equation}
	\end{definition}	
	\noindent The two tensors   $\mathcal {X},\mathcal {Y}\in \mathbb{R}^{I_{1}\times I_{2}\times \ldots \times I_{N}\times   J_{1}\times J_{2}\times \ldots \times J_{M}}$ are orthogonal iff 
	$	\langle \mathcal {X},\mathcal {Y} \rangle=0$.\\
	
	\begin{proposition}
		Let $\mathcal {A} \in \mathbb{R}^{I_{1}\times \ldots \times I_{N}\times   K_{1}\times \ldots \times K_{N}}$, 
		$\mathcal {B} \in \mathbb{R}^{K_{1}\times \ldots \times K_{N}\times   J_{1}\times \ldots \times J_{M}}$. Then
		\begin{enumerate}
			\item $(\mathcal {A}\ast_N B)^T= \mathcal {B}^T  \ast_N  \mathcal {A}^T$.
			\item $\mathcal {I_N}  \ast_N \mathcal {B}=\mathcal {B}$ and $\mathcal {B}  \ast_M \mathcal {I_M}=\mathcal {B}$, where identity tensors $\mathcal {I_N} $ and $\mathcal {I_M} $ are such that 
			$\mathcal {I_N} \in \mathbb{R}^{K_{1} \times\ldots \times K_{N}\times   K_{1}\times \ldots \times K_{N}}$ and 
			$\mathcal {I_M} \in \mathbb{R}^{J_{1}\times \ldots \times J_{M}\times   J_{1} \times \ldots \times J_{N}}$.
			
		\end{enumerate}
		
	\end{proposition}
	
	\medskip
	\section{Tensor formulation using RBF descritization}
	Consider linear steady problem
	\begin{equation}\label{steady}
	\mathcal{L} u(x)=f(x),\quad x\in\Omega\subset\mathbb{R}^3,
	\end{equation}
	where $\mathcal{L}$ is a linear differential operator. Equation (\ref{steady}) is subject to a homogeneous condition on its boundary
	$\partial \Omega$ of the form 
	\begin{equation}\label{boundary}
	\mathcal{ B}u(x)=g,\quad x\in\partial\Omega,
	\end{equation}
	The most common choice of MQ-RBF is given as the following
	\begin{equation}
	\varphi_\varepsilon(r)=\sqrt{1+\varepsilon^2r^2},
	\end{equation}
	where $r=\|x\|_2$, $x\in\mathbf{R}^3$, is the argument that makes $\phi$ radially symmetric about its center  and $\varepsilon$ is referred to as the shape parameter. To ensure the existence of a solution to the
	boundary problem (\ref{steady}), the boundary $\partial\Omega$ is supposed to be
	sufficiently smooth. An appropriate choice of the shape  parameter $\varepsilon$ makes MQ the best candidate for a good approximation among all the other RBF choices; see \cite{SK}. Given the collocation points $\left\{\left(x_{m n p}, y_{m np},z_{mnp}\right)\right\}_{m=1, n=1,p=1}^{M, N,P}$, the MQ collocation method suggests for each point $\left(x_{i jk}, y_{i jk},z_{ijk}\right), i=1, \ldots, M, j=$ $1, \ldots, N,$, $k=1, \ldots, P$, the following approximant
	\begin{equation}\label{approximat}
	u_\text{app}\left(x_{i jk}, y_{i jk},z_{ijk}\right)=\sum_{m=1}^{M} \sum_{n=1}^{N} \sum_{p=1}^{P} \alpha_{mnp} \varphi\left(r_{m n p}^{ijk}\right), 
	\end{equation}
	for $ m=1, \ldots, M, n=$ $1, \ldots, N,$, $p=1, \ldots, P$, and
	\begin{equation}
	\varphi_\varepsilon\left(r_{m n}^{i jk}\right)=\varphi_\varepsilon\left(\left\|\left(x_{i jk}, y_{i jk},z_{ijk}\right)-\left(x_{m n p}, y_{m np},z_{mnp}\right)\right\|\right).
	\end{equation}
	The expansion coefficients $\alpha_{mnp}$ are determined by enforcing the interpolation condition
	\begin{equation}
	u_\text{app}\left(x_{m n p}, y_{m np},z_{mnp}\right)=u\left(x_{m n p}, y_{m np},z_{mnp}\right).
	\end{equation}
	By using Definition \ref{einstein}, a tensor equation with Einstein product associated
	with (\ref{approximat}, can be defined as follows,
	\begin{equation}\label{basis}
	\mathcal{A}\ast_3\Gamma=\mathcal{U},
	\end{equation}
	where $\mathcal{A}\in\mathbb{R}^{M\times  N\times P\times M\times N\times P}$ is a sixth order tensor with entries
	\begin{equation}
	\left(\mathcal{A}\right)_{ijkmnp}=\varphi\left(r_{m n p}^{ijk}\right),  m,i=1, \ldots, M; n,j=1, \ldots, N; p,k=1, \ldots, P,
	\end{equation}
	and 
	\begin{equation}
	(\Gamma)_{mnp}=\alpha_{mnp},\quad (\mathcal{U})_{mnp}=u\left(x_{m n p}, y_{m np},z_{mnp}\right),
	\end{equation}
	for $ m=1, \ldots, M, n=$ $1, \ldots, N,$ $p=1, \ldots, P$. We refer to $\mathcal{A}$ as the system tensor  serving as the basis of the approximation space.\\
	
	\noindent Now, let $\mathcal{B}$ be the linear operator associated with the boundary conditions (\ref{boundary}),  and consider that the collocation points  $\left\{\left(x_{m n p}, y_{m np},z_{mnp}\right)\right\}_{m=1, n=1,p=1}^{M, N,P}$ subsets.
	One subset contains $\left\{\left(x_{m n p}, y_{m np},z_{mnp}\right)\right\}_{m=1, n=1,p=1}^{M_\mathcal{L}, N_\mathcal{L},P_\mathcal{L}}$, where the PDE is enforced and the other
	subset $\left\{\left(x_{m n p}, y_{m np},z_{mnp}\right)\right\}_{m=1, n=1,p=1}^{M_\mathcal{B}, N_\mathcal{B},P_\mathcal{B}}$, where boundary conditions are enforced.\\
	In the MQ collocation method and when applying the linear operator $\mathcal{L}$, we take for each point $\left(x_{ijk}, y_{ijk},z_{ijk}\right)$, for $ i=1, \ldots, M_{\mathcal{L}}, j=1, \ldots, N_{\mathcal{L}}$, $k=1, \ldots, P_{\mathcal{L}}$
	\begin{equation}\label{operL}
	\mathcal{L}u\left(x_{ijk}, y_{ijk},z_{ijk}\right)=\sum_{m=1}^{M} \sum_{n=1}^{N} \sum_{p=1}^{P} \alpha_{mnp} \mathcal{L}\varphi\left(r_{mnp}^{ijk}\right), 
	\end{equation}
	and when applying the operator $\mathcal{B}$ for $ i=M_{\mathcal{L}}+1, \ldots, M, j=N_{\mathcal{L}}+1, \ldots, N$, $k=P_{\mathcal{L}}+1, \ldots, P$
	\begin{equation}\label{operB}
	\mathcal{B}u\left(x_{ijk}, y_{ijk},z_{ijk}\right)=\sum_{m=1}^{M} \sum_{n=1}^{N} \sum_{p=1}^{P} \alpha_{mnp} \mathcal{B}\varphi\left(r_{m n p}^{ijk}\right).
	\end{equation}
	In tensor notation (Definition \ref{einstein}), the right hand-side of equations
	(\ref{operL}) and (\ref{operB})) can be written as $\mathcal{H}\ast_{N}\Gamma$, where for $ m=1, \ldots, M, n=$ $1, \ldots, N,$ $p=1, \ldots, P$, the entries of the tensor $\mathcal{H}$ are defined as follows
	\begin{equation*}
	\left(\mathcal{H}\right)_{ijkmnp}=\mathcal{L}\varphi\left(r_{m n p}^{ijk}\right),  \quad i=1, \ldots, M_{\mathcal{L}}, j=1, \ldots, N_{\mathcal{L}}, k=1, \ldots, P_{\mathcal{L}},
	\end{equation*}
	and
	\begin{equation*}
	\left(\mathcal{H}\right)_{ijkmnp}=\mathcal{L}\varphi\left(r_{m n p}^{ijk}\right),\quad i=M_{\mathcal{L}}+1, \ldots, M, j=N_{\mathcal{L}}+1, \ldots, N, k=P_{\mathcal{L}}+1, \ldots, P.
	\end{equation*}
	
	\noindent By inverting the system tensor in (\ref{basis}), the tensor
	that discretizes the PDE in space is the differentiation tensor
	\begin{equation}
	\mathcal{T}=\mathcal{H}\ast_3\mathcal{A}^{-1}.
	\end{equation}
	The linear  steady problem (\ref{steady}) with boundary condition (\ref{boundary}) is discretized as
	\begin{equation}\label{tensoreq}
	\mathcal{T}\ast_3\mathcal{U}=\mathcal{F},
	\end{equation}
	where $(\mathcal{F})_{ijk}=f\left(x_{ijk}, y_{ijk},z_{ijk}\right),$ $i=1, \ldots, M, j=1, \ldots, N, k=1, \ldots, P.$
	The problem (\ref{steady})-(\ref{boundary}) has a solution
	\begin{equation}
	\mathcal{U}=\mathcal{ A}\ast_3\mathcal{ Y},
	\end{equation}
	where $\mathcal{ Y}$ is the solution of the following ill-posed  tensor equation
	\begin{equation}\label{tensorprob}
	\mathcal{ H}\ast_3\mathcal{ Y}=\mathcal{F}.
	\end{equation}
	
	\section{Hierarchical MQ-RBF Interpolation}
	It can be easily seen from the MQ-RBF formulation, that the tensors formulated 
	are fully populated. This therefore leads to  tensors $\mathcal{ A}$ in (\ref{basis}) and $\mathcal{ H}$ in (\ref{tensorprob}) that cannot be stored in memory. Limitation in memory storage 
	makes the MQ-RBF approach unattractive and limits its applicability on classical computers when the number of degrees of freedom reaches a few thousand,
	which is often not sufficient in practice.  For the matrix case, many applications such as the matrices computed for the boundary element method (BEM), 
	fast convolution techniques  with usual  kernels on some unstructured grids,  have been developed to perform   the matrix-vector products in a
	reasonable time leading to fast solvers. It is also our aim to develop fast linear solvers for the tensor case using the MQ-RBF formulation. Based on algebraic compression proposed in  \cite{Aussal}, we will be able to provide new fast linear solvers that require much less memory storage for solving  large linear tensor equations obtained in the MQ-RBF framework. To this aim, we extend for the tensor case, the method of algebraic compression based on divide and
	conquer process introduced for the matrix case to accurately approximate the full matrix with a hierarchical one with low
	rank pieces \cite{Aussal}. Therefore, in a similar manner to the matrix case, we follow three steps to build the compressed tensors in equations (\ref{basis}) and  (\ref{tensorprob}).
	
	For the first step, we use a binary domain decomposition to
	compute two independent binary trees  for the three-dimensional set represented by the collocation points. To
	keep a well balanced spatial distribution with any spatial configuration, 
	geometric and median cutting approaches are used, which provide the
	best way for all the groups of collocation points  encountered at
	each depth of the tree. The subdivision of the collocation points is carried out until the number of points in a group falls below the low a threshold  $T_{\text{leaf}}\left(\log (M)+\log (N)+\log (P)\right)^{\frac{3}{2}}.$\\
	At the second step, the binary tree domains  associated
	to the collocation points, allow block interactions for algebraic compression. To proceed the hierarchical construction of our compressed tensors, the
	blocks are defined by the sets of collocation points  $X_I$ defined by the particles $\left(x_{i jk}, y_{i jk},z_{ijk}\right)$ for $i,j,k\in I$ and $Y_J$ defined by the particles $\left(x_{m n p}, y_{m np},z_{mnp}\right)$
	for $m,n,p\in J$. For these two sets of collocation points, the entries of the tensor $\mathcal{ A}$ in equation (\ref{basis}) are of the form
	\begin{equation}
	\left(\mathcal{A}\right)_{ijkmnp}=\varphi_\varepsilon\left(\left\|\left(x_{i jk}, y_{i jk},z_{ijk}\right)-\left(x_{m n p}, y_{m np},z_{mnp}\right)\right\|\right),  i,j,k\in I, m,n,p\in J.
	\end{equation}
	Moreover, we choose  two  axis-parallel boxes $B_I$ and $B_J$ that surrounding each set of particles $X_I$ and $Y_J$, 
	respectively. If the bounding boxes satisfy
	the following admissibility condition
	$$\max \left\{\operatorname{diam}\left(B_{I}\right), \operatorname{diam}\left(B_{J}\right)\right\} \leq \eta \operatorname{dist}\left(B_{I}, B_{J}\right),$$
	for fixed $\eta>0,$ then the MQ-RBF admits the degenerated expansion
	\begin{eqnarray}
	&&\varphi_\varepsilon\left(\left\|\left(x_{i jk}, y_{i jk},z_{ijk}\right)-\left(x_{m n p}, y_{m np},z_{mnp}\right)\right\|\right) \approx\\ &&\sum_{\mu=1}^{p} \sum_{\nu=1}^{p} L_{r, \mu}(x) g\left(\xi_{r, \mu}, \xi_{c, \nu}\right) L_{c, \nu}(y), \quad x \in X_{r}, y \in Y_{c}
	\end{eqnarray}
	At the third and the last step, the Adaptive Cross Approximation
		\cite{Liu2020} is used for admissible interactions, completed by the standard full computation for close interactions. For the evaluation of the convergence of the ACA algorithm, we  use the  similar  criterion in \cite{Liu2020}. The distances between each set of particles ${X}_{I}$ and ${Y}_{J}$ are evaluated from their projections on the axis defined by the two centres of each dataset.
	
	\section{Tensor Krylov subspace methods}
	In this section, we propose iterative methods based on Global Arnoldi and Global Golub– Kahan bidiagonlization (GGKB), combined with Tikhonov regularization, to solve the ill-posed tensor equation (\ref{tensoreq}).   Solving (\ref{tensoreq}) is equivalent to finding the solution of the following minimization problem
	\begin{equation}\label{lsq}
	\min_{\mathcal{Y}}\|\mathcal{H}\ast_3\mathcal{Y}-\mathcal{F}\|_F.
	\end{equation}
	In order to diminish the effect of the ill-conditioning of the tensor $\mathcal{H}$, we replace the original problem  by a better stabilized one. One of the most popular regularization
	methods is due to Tikhonov \cite{tikhonov}. The method replaces the problem  by the new one
	\begin{equation}\label{Tikh}
	\min_{\mathcal {Y} } \Vert \mathcal {H} \ast_3  \mathcal {Y} -  \mathcal {F} \Vert_F^2 + \mu \Vert \mathcal{D}  \ast_N  \mathcal {Y} \Vert_F^2,
	\end{equation}
	where $\mu\geq0$ is the regularization parameter and $\mathcal{D}$ is a regularization tensor chosen to obtain a solution with desirable properties. The tensor $\mathcal{ D}$ could be the identity tensor or a discrete form of first or second derivative. In the first case, the parameter $\mu$ acts on the size of the solution, while in the second case $\mu$ acts on the smoothness of the solution. We will only consider  the particular case where the tensor $\mathcal{D}$ reduces to the identity tensor $\mathcal{I}_N$. Therefore, Tikhonov regularization problem in this case is of the following form
	\begin{equation}\label{Tikhrd}
	\mathcal{Y}_\mu=\text{arg}\underset{ \mathcal {Y}}{ \text{min}} \left(\| \mathcal {H} \ast_3  \mathcal {Y}  -\mathcal F\|_{F}^2+\mu\|\mathcal{Y}\|_F^2 \right),
	\end{equation}
	Many techniques for choosing a suitable value of $\mu$ have been analyzed and illustrated in the literature; see, e.g., \cite{wahbagolub} and references therein. In this paper we will use the discrepancy principle and the Generalized Cross Validation (GCV) techniques.
	
	\subsection{Tensor Global GMRES method}
	Let  $\mathcal {V} \in \mathbb{R}^{M\times N\times P}$ and consider the following $m$-th tensor Krylov subspace  defined as
	\begin{equation}\label{krylov1}
	\mathcal {K}_m(\mathcal {H},\mathcal {V})=span\{\mathcal {V},\mathcal {H}\ast_3\mathcal {V},\ldots,\mathcal {H}^{m-1}\ast_3\mathcal {V} \},
	\end{equation}
	where $\mathcal {H}^i\ast_3\mathcal {V}=\mathcal {H}\ast_3 \mathcal {H}^{i-1}\ast_3\mathcal {V}$. The tensor global Arnoldi algorithm can easily defined as the global Arnoldi process for defined in \cite{jbilou1} for the matrix case. The algorithm is defined as follows (see \cite{beik1,Elguide,huang1,jbilou1})
	
	\begin{algorithm}
		\caption{Einstein Tensor  Global Arnoldi process (ETGA)}\label{alg1}
		\begin{enumerate}
			
			\item Inputs: A tensor $\mathcal {H} \in\mathbb{R}^{M\times N\times P \times M\times N\times P }$, and a tensor   $\mathcal {V} \in \mathbb{R}^{M\times N\times P}$ and the integer $m$.
			\item Set $\beta=\Vert \mathcal {V} \Vert_F$ and $\mathcal {V}_1=\mathcal {V} /\beta$.
			\item For $j=1,\ldots,m$
			\item $\mathcal {W}=\mathcal {H} \ast_3\mathcal {V}_j$
			\item for $ i=1,\ldots,j$.
			\begin{itemize}
				\item $h_{i,j}=\langle \mathcal {V}_i,\mathcal {W} \rangle$,
				\item $\mathcal {W}= \mathcal {W}-h_{i,j}\mathcal {V}_i$
			\end{itemize}
			\item endfor
			\item $h_{j+1,j}=\Vert \mathcal {W} \Vert_F$. If $h_{j+1,j}=0$, stop; else
			\item $\mathcal {V}_{j+1}=\mathcal {W}/h_{j+1,j}$.
			\item EndFor
			
		\end{enumerate}
	\end{algorithm}
	
	\medskip 
	\noindent 
	Application of $m$  steps of Algorithm \ref{alg1}
	yields the decompositions
	\begin{equation}
	\label{lin1}
	\mathcal {H} \ast_3 \mathbb{V}_m = \mathbb{V}_{m+1} \times_{4} {\widetilde H}_m^T,
	\end{equation}
	where $\mathbb{V}_m$ is the 4-mode tensor with frontal slices $\mathcal {V}_1,\mathcal {V}_2,\ldots,\mathcal {V}_m$ and   $\mathcal {H} \ast_3 \mathbb{V}_m$ is the $4$-mode tensor with frontal slices $\mathcal {H} \ast_N \mathcal {V}_1,\ldots,\mathcal {H} \ast_N \mathcal {V}_m$ and $\widetilde{H}_m \in \mathbf{R}^{(m+1) \times m}$ is the following upper Hessenberg matrix 
	$$\widetilde{H}_{m}=\left[\begin{array}{ccccc}
	h_{1,1} & h_{1,2} & h_{1,3} & \dots & h_{1, m} \\
	h_{2,1} & h_{2,2} & h_{2,3} & \dots & h_{2, m} \\
	0 & h_{3,2} & h_{3,3} & \dots & h_{3, m} \\
	\vdots & \ddots & \ddots & \ddots & \vdots \\
	\vdots & & 0 & h_{m, m-1} & h_{m, m} \\
	0 & \dots & \dots & 0 & h_{m+1, m}
	\end{array}\right].$$ 
	\noindent Notice  that the
	$\mathcal {V}_i$'s obtained from Algorithm \ref{alg1} form an orthonormal basis of the tensor Krylov subspace $\mathcal {K}_m(\mathcal {H},\mathcal {V})$. 
	We can now define the Tensor GMRES method to solve the problem \eqref{Tikhrd}. Using  the global GMRES, we look for an approximate solution $\mathcal {X}_m$, starting from $\mathcal {X}_0$ such that  $\mathcal {X}_m \in \mathcal {X}_0+ \mathcal {K}_m(\mathcal {H},\mathcal {R}_0)$, with $\mathcal{R}_0=\mathcal{H}\ast_3\mathcal{X}_0-\mathcal{ F}$. 
	\noindent Using the relation \eqref{lin1}, we can show that (see \cite{Elguide,huang1})
	\begin{equation}\label{tgmres2}
	\mathcal {X}_m=\mathcal {X}_0+ \mathbb{V} \bar{\times}_{4} y_m,
	\end{equation}
	where $y_m\in\mathbb{R}^{m}$. Therefore, replacing (\ref{tgmres2}) into (\ref{Tikhrd}), yields the reduced minimization problem 
	
	\begin{equation}\label{tgmres3}
	y_{m,\mu}=\arg \min_y \left(\Vert \widetilde H_m y- \beta e_1\Vert_2^2+\mu\Vert y\Vert_2^2\right),
	\end{equation}
	where $\beta=\left\|\mathcal{R}_{0}\right\|_{F}$   and $e_{1}=(1,0, \ldots, 0)^{T} \in \mathbb{R}^{m+1}$. The solution \eqref{tgmres3} can obtained by the solving the following reduced least square problem is given by
	
	\begin{equation}
	\label{min}
	y_{m,\mu}= \arg \min_y\left  \Vert \left ( \begin{array}{ll}
	\widetilde H_m\\
	\mu I_m
	\end{array}\right ) y - \left ( \begin{array}{ll}
	\beta e_1\\
	0
	\end{array}\right ) 
	\right \Vert_2.
	\end{equation}
	The minimizer $y_{m,\mu}$  of the problem \eqref{min} is computed as the solution of the linear system of equations
	\begin{equation}{\label{min1}}
	\widetilde H_{m,\mu} y=\widetilde H_m^T\beta e_1,
	\end{equation}
	where $\widetilde H_{m,\mu}= (\widetilde H_m^T \widetilde H_m+ \mu^2 I_m)$.

	\noindent The Tikhonov problem \eqref{tgmres3} is a matrix one with small dimension as $m$ is generally small and then can be solved by some techniques such as the GCV method \cite{golubwahba} or the L-curve criterion \cite{hansen1,hansen2,reichel1,reichel2}. \\
	An appropriate selection of the regularization parameter $\mu$
	is important in  Tikhonov regularization. Here we can use  the
	generalized cross-validation (GCV) method
	\cite{golubwahba,wahbagolub}. For this method, the regularization
	parameter is chosen to minimize the GCV function
	$$GCV(\mu)=\frac{||\widetilde H_m y_{m,\mu}-{\bf
			\beta e_1}||_2^2}{[tr(I_m-\widetilde H_m  \widetilde H_{m,\mu}^{-1}\widetilde H_m^T)]^2}=\frac{||(I_m-\widetilde H_m \widetilde H_{m,\mu}^{-1} \widetilde H_m^T){\beta e_1}||_2^2}{[tr(I_m-H_m H_{m,\mu}^{-1} \widetilde H_m^T)]^2}.$$  As the projected problem we are dealing with is of small size, we cane use the SVD decomposition of $\widetilde H_m$ to obtain a more simple and computable expression of $GCV(\mu)$. Consider the SVD decomposition of $\widetilde H_m=U\Sigma V^T$. Then the GCV function could be expressed as (see \cite{wahbagolub})
	\begin{equation}
	\label{gcv2}
	GCV(\mu)=\frac{\displaystyle
		\sum_{i=1}^m\left(\frac{\tilde
			g_i}{\sigma_i^2+\mu^2}\right)^2}{\displaystyle\left(\sum_{i=1}^m
		\frac{1}{\sigma_i^2+\mu^2}\right)^2},
	\end{equation}
	where $\sigma_i$ is the $i$th singular value of the matrix
	$\widetilde H_m$ and $\tilde g= \beta U^T e_1$.
	\\
	In  the practical implementation, it's more convenient to use a restarted version of the Global GMRES. 
	The tensor Global GMRES  algorithm for solving tensor linear  equations (\ref{tensorprob}) is summarized as
	follows:
	\vspace{0.3cm}
	\begin{algorithm}[!h] 
		\caption{ Einstein Tensor Global GMRES method for Tikhonov regularization}\label{TG-GMRES}
		\begin{enumerate}
			\item {\bf Inputs:}  tensors $\mathcal {H}$, $\mathcal {F}$, initial guess $\mathcal{U}_0$, a tolerance $tol$, number of iterations between restarts $m$ and \textbf{Maxit}: maximum number of outer iterations.
			\item Compute $\mathcal{R}_0=\mathcal {H} \ast_3 \mathcal {U}_0  -\mathcal F$, set $\mathcal{V}=\mathcal{R}_0$ and $k=0$.
			\item Determine the orthonormal bases  $\mathbb{V}_{m}$ of tensors, and the upper Hessenberg matrix  $\widetilde{H}_m$
			matrix by applying Algorithm \ref{alg1} to the pair $\left(\mathcal{H}, \mathcal{ R}_0\right)$ 
			\item  Determine $\mu_{k}$  as the parameter minimizing the GCV function  given by (\ref{gcv2})
			\item  Determine $y_m$ as a solution of the reduced Tikhonov regularization problem (\ref{tgmres3}) and then compute $\mathcal{X}_m$ by (\ref{tgmres2}) \item  If $\left\| \mathcal{R}_0\right\|_{F}<tol$ or $k>\textbf{Maxit}$; Stop \\
			else: set $\mathcal X_{0}=\mathcal{X}_{m}$, $k=k+1,$ goto 2
		\end{enumerate}
	\end{algorithm}
	
	\subsection{Einstein Tensor Global LSQR}
	Instead of using the tensor global Arnoldi to generate a basis for the projection subspace, we can use the tensor global Golub-Kahan  process.   Given $\mathcal{ H}$  and $\mathcal{F}$, the tensor global Golub-Kahan  algorithm is defined as follows
	
	\begin{algorithm}[!h]\label{alg3}
		\caption{ Einstein Tensor Global Golub Kahan algorithm}
		\begin{enumerate}
			\item {\bf Inputs} The tensors $\mathcal {H}$, $\mathcal {F}$,  and an integer $m$.
			\item 	Set $\sigma_1= \Vert \mathcal{F} \Vert_F$, 
			$\mathcal {P}_1=\mathcal {F}/\sigma_1$ and  $\mathcal {Q}_1=0$
			\item for $j=2,\ldots,m$
			\begin{enumerate}
				\item $\widetilde {\mathcal {Q}}= \mathcal {H} \ast_3 \mathcal {P}_{j-1} -\sigma_{j-1}\mathcal {Q}_{j-1}$
				\item $\rho_j=\Vert \widetilde {\mathcal {Q}}\Vert_F$ if $\rho_j=0$ stop, else
				\item $\mathcal {Q}_j=\widetilde {\mathcal {Q}}/\rho_j$
				\item $\widetilde {\mathcal {P}}=\mathcal {H}^T \ast_3 \mathcal {Q}_j-\rho_j \mathcal {P} _{j-1}$
				\item $\sigma_j=\Vert \widetilde {\mathcal {P}} \Vert_F$
				\item if $\rho_j=0$ stop, else
				\item $\mathcal {P}_j=\widetilde {\mathcal {P}}/\sigma_j$
			\end{enumerate}
		\end{enumerate}
		
	\end{algorithm}
	
	\medskip 
	
	\noindent 
	Application of  $m$ steps of the GGKB method to $\mathcal{H}$ with
	initial tensor $\mathcal{F}$,  produces the lower bidiagonal matrix $C_m\in\mathbf{R}^{m\times m}$ 
	\[
	C_m=\begin{bmatrix}
	\rho_1\\
	\sigma_2 & \rho_2&\\
	&\ddots&\ddots\\
	& & \sigma_{m-1}&\rho_{m-1}\\
	&&&\sigma_m&\rho_m
	\end{bmatrix}
	\]
	and 
	\[
	\widetilde{C}_m=\begin{bmatrix}
	C_m\\
	\sigma_{m+1}e_m^T
	\end{bmatrix}\in\mathbb{R}^{(m+1)\times m}.
	\]
	If we assume that $m$ is small enough so that all nontrivial entries of the matrix $\widetilde{C}_{m}$ are positive, then Algorithm \ref{alg3} yield   the decompositions
	\begin{eqnarray}
	\mathcal {H} \ast_3 \mathbb{Q}_m& = &\mathbb{P}_{m+1} \times_{4} {\widetilde C}_m,\\
	\mathcal {H }^T \ast_3 \mathbb{P}_m	& = &\mathbb{Q}_{m} \times_{4} { C}_m^T,
	\end{eqnarray} where $\mathbb{Q}_m$  and $\mathbb{P}_{m+1}$ are  4-mode tensors with orthonormal frontal slices $\mathcal {Q}_1,\ldots,\mathcal {Q}_m$ \\and  $\mathcal {P}_1,\mathcal {P}_2,\ldots,\mathcal {P}_{m+1}$, respectively.
	\noindent 
	In global LSQR, the approximate solution is defined as
	\begin{equation}\label{Xkmu}
	\mathcal {X}_m=\mathbb{Q}_m \bar{\times}_{4} y_m,
	\end{equation}
	where $y_{m}$ solves 
	\begin{equation}\label{normeq2}
	(\widetilde{C}_m^T\widetilde{C}_m+\mu I_m)y=\sigma_1\widetilde{C}_m^Te_1,\qquad\sigma_1=\|\mathcal{F}\|_F.
	\end{equation}
	It is also computed by solving the least-squares problem
	\begin{equation}\label{leastsq}
	\min_{y\in\mathbb{R}^m} \begin{Vmatrix}
	\begin{bmatrix}
	\mu^{1/2}\widetilde{C}_m\\
	I_m
	\end{bmatrix}
	y-\sigma_1\mu^{1/2}e_1 \end{Vmatrix}_2.
	\end{equation}
	The following algorithm summarizes the main steps to compute a regularization parameter and the corresponding regularized solution of (\ref{Tikhrd}), using Einstein Tensor GGKB for Tikhonov regularization.
	
	\begin{algorithm}[!h]
		\caption{ Einstein Tensor global LSQR method for Tikhonov regularization}\label{TG-GKB}
		\begin{enumerate}
			\item {\bf Inputs}: The tensors $\mathcal {H}$, $\mathcal {F}$ and $m$.
			\item Determine the orthonormal bases $\mathbb{P}_{m+1}$ and $\mathbb{Q}_{m}$ of tensors, and the bidiagonal $C_m$ and $\widetilde{C}_m$
			matrices  with Algorithm\ref{alg3}.
			\item Determine $\mu$ that minimize the GCV function.
			\item  Determine $y_{m,\mu}$ by solving  (\ref{leastsq}) and then compute $\mathcal{ X}_{m,\mu}$ by (\ref{Xkmu}).
		\end{enumerate}
		
	\end{algorithm}
	\section{Numerical results}
	This section performs some numerical tests on the methods of Tensor Global GMRES(m) and Tensor Global LSQR algorithm given by Algorithm \ref{TG-GMRES} and Algorithm \ref{TG-GKB}, respectively, 
	for solving tensor equation in the form
	(\ref{tensorprob}) resulting from RBF disretization of (\ref{steady})-(\ref{boundary}). We evaluate our proposed methods to solve the three-dimensional (3D) acoustic Helmholtz equation
	\begin{equation}\label{Helm}
	\Delta u(x)+k^{2} u(x)=0, \quad x \in \Omega \subset \mathbb{R}^{3},
	\end{equation}
	where $\Delta$ is the Laplace operator, $u(x)$ is the sound pressure at point $x, k=\omega / c$ is the wave number with the circular frequency $\omega$ and the speed of sound $c$ through the fluid medium. Equation (\ref{Helm}) is subject to a homogeneous condition on its boundary an of the form
	\begin{equation}\label{Helmboun}
	a(x) u(x)+b(x) \frac{\partial u(x)}{\partial n}=0, \quad x \in \partial\Omega,
	\end{equation}
	where $\displaystyle \frac{\partial}{\partial n}$ denotes the outward normal to the boundary at point $x$. It's known that different
	distribution collocation points affect the results of the RBF-Based Meshless Method. In order to
	measure the efficiency and accuracy of the proposed algorithms, different types of collocations point distributions
	are considered. The sets of collocation points we considered are random, uniform and Halton. For the results reported in the examples,  we used the stopping criterion given by,
	\[
	\frac{\left\|\mathcal{X}_{m+1,\mu}-\mathcal{X}_{m,\mu}\right\|_F}{\left\|\mathcal{X}_{m,\mu}\right\|_F} \leq \tau.
	\]
	The maximum number of 200 iterations was allowed for both algorithms. To determine the effectiveness of our solution methods, we evaluate 
		$$\text{Relative error}=\frac{\left\|{ \mathcal U}-{\mathcal U}_{\text{computed}}\right\|_{F}}{\|{ \mathcal U}\|_{F}}.$$ 
	All computations were carried out using the MATLAB environment on an Intel(R) Core(TM) i7-8550U CPU @ 1.80GHz (8 CPUs) computer with 12 GB of
	RAM. The computations were done with approximately 15 decimal digits of relative
	accuracy. 
	\subsection{Example 1}
\begin{figure}
		\begin{center}
			\includegraphics[width=.33\textwidth,height=.35\textwidth]{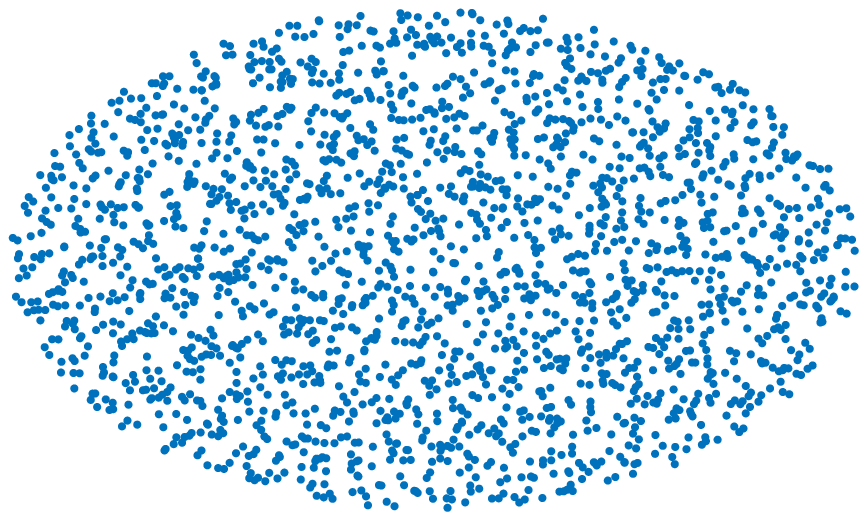}%
			\includegraphics[width=.33\textwidth,height=.33\textwidth]{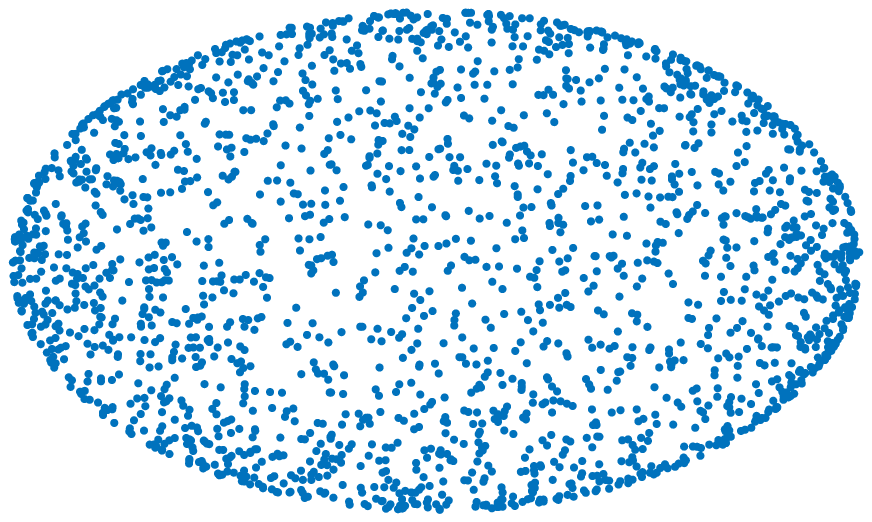}%
			\includegraphics[width=.33\textwidth,height=.34\textwidth]{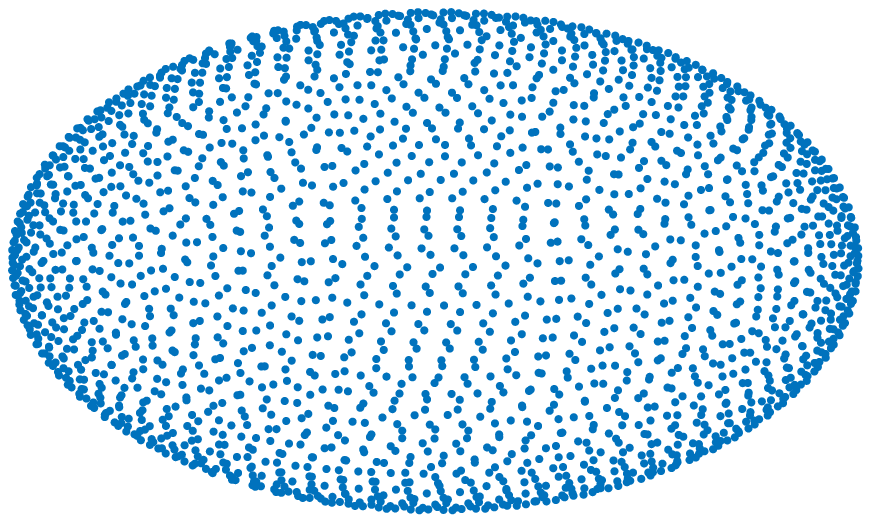}
			\caption{Sets of collocation points:  random (left),  uniform (middle) and  Halton (right)} 
			\label{fig1}
		\end{center}
	\end{figure}

\begin{figure}
		\begin{center}
			\includegraphics[width=.5\textwidth,height=.4\textwidth]{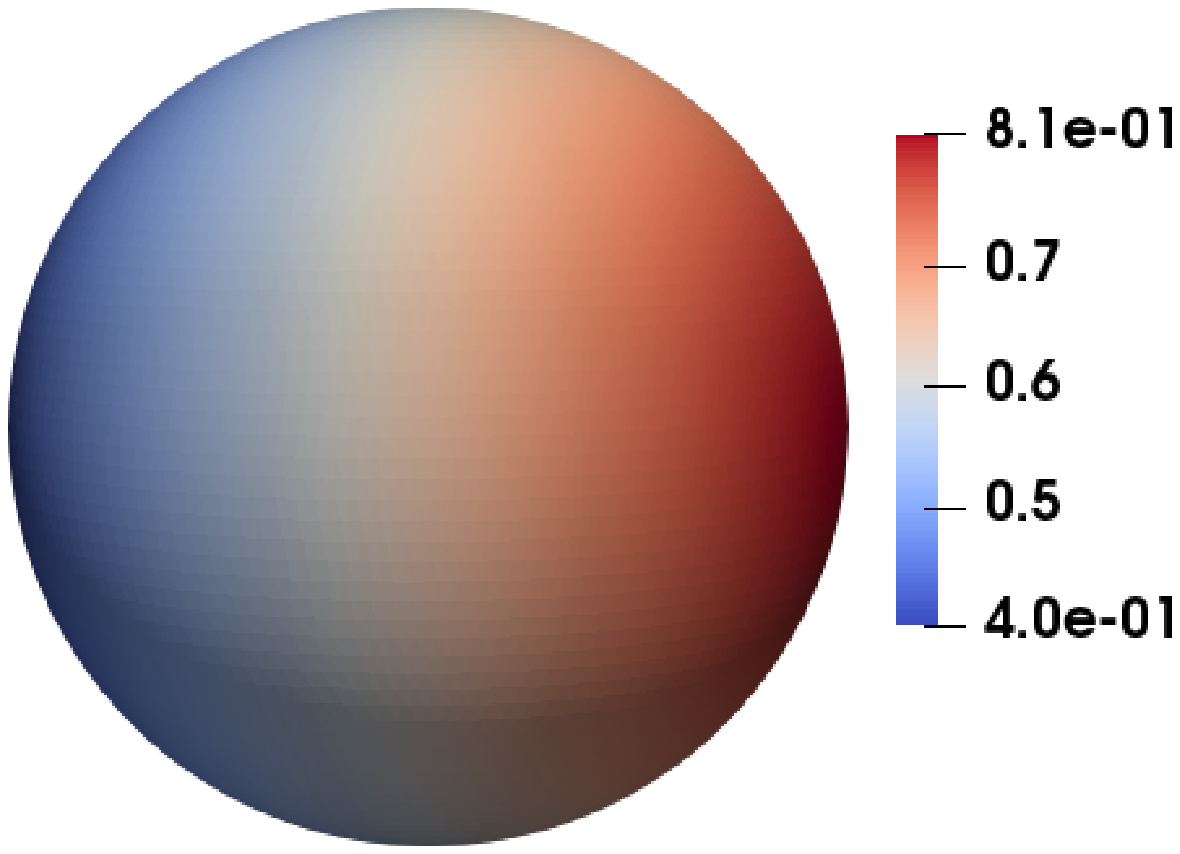}%
			\includegraphics[width=.5\textwidth,height=.4\textwidth]{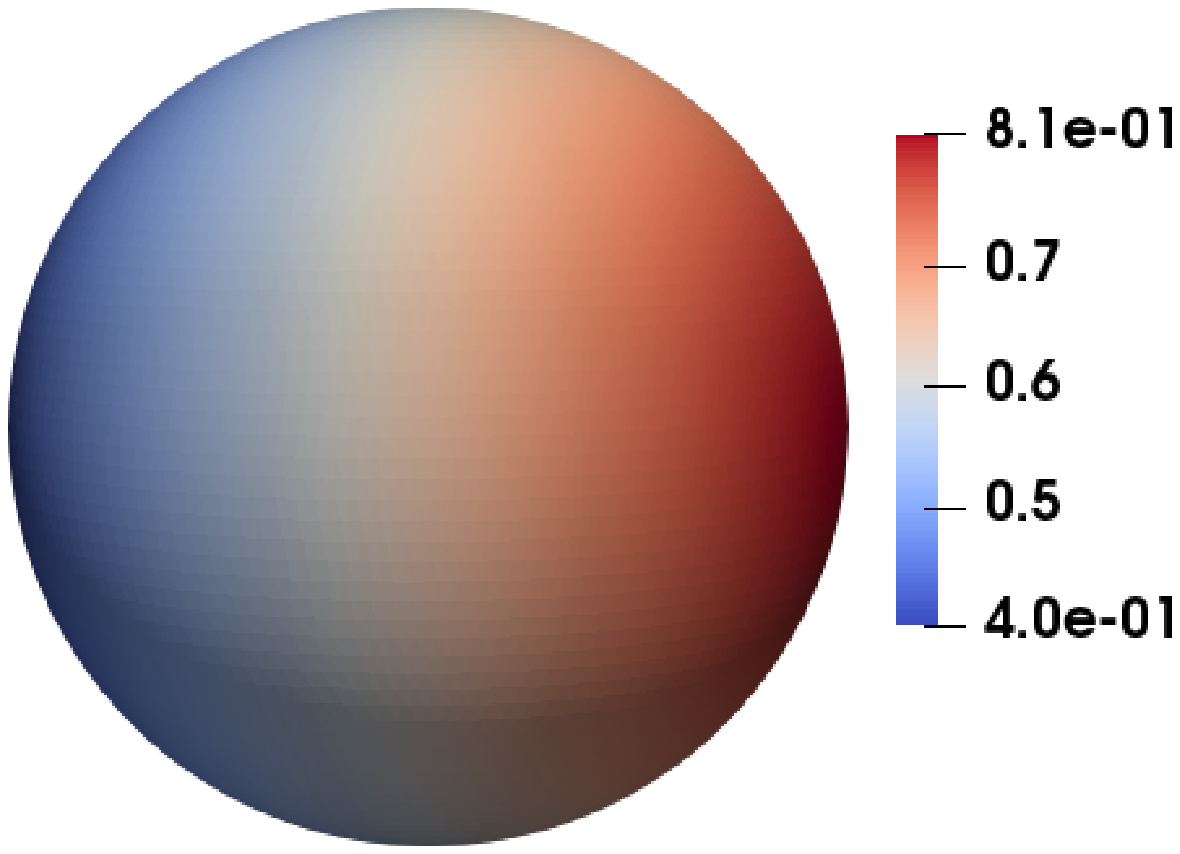}
			\caption{Exact solution (left) and approximate solution (right).} \label{fig2}
		\end{center}
	\end{figure}

\begin{figure}
		\begin{center}
			\includegraphics[width=.8\textwidth,height=.4\textwidth]{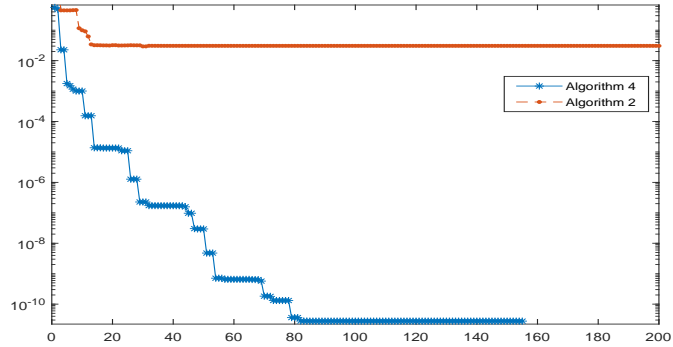}%
			\caption{Example 1: Values of the relative error 
				versus the number of iterations.} \label{REex1}
		\end{center}
	\end{figure}

	In this example we are interested in the numerical solution of  three-dimensional Helmholtz problem  (\ref{Helm}) in a unit
	sphere domain  with homogeneous Dirichlet boundary conditions,
	i.e., (\ref{Helmboun}) with $a(x)=1$ and $b(x)=0 ,$ using Algorithm \ref{TG-GMRES}  and Algorithm \ref{TG-GKB}. The three different distributions of the collocation points are shown in Figure \ref{fig1}. The chosen exact solution is the following
	\begin{equation}
	u(x,y,z)=\text{exp}\left(-\frac{(x-0.25)^2+(y-0.25)^2+z^2}{\sigma}\right),
	\end{equation}
	with $\sigma=20$. The exact solution of the problem is reported on the left of Figure \ref{fig2}. In Table \ref{tab1}, we evaluate the effectiveness of
	Algorithm \ref{TG-GMRES} and Algorithm \ref{TG-GKB}, when solving the three-dimensional Helmholtz problem for different
	numbers of distribution evaluation points. We denote by "Iter"  the iteration steps and it is obtained as soon as $\tau$ reaches $10^{-12}$. In Figure \ref{REex1}, we plot the values of the relative error  versus the number of iterations, obtained by applying Algorithm \ref{TG-GMRES} and Algorithm \ref{TG-GKB} for the random distribution of collocation points with $N=M=P=12$.
	\begin{table}[htbp]
		\begin{center}\footnotesize
			\renewcommand{\arraystretch}{1.3}
			\begin{tabular}{cccccc}\hline
				\multicolumn{1}{c}{{Collocation points}} &\multicolumn{1}{c}{{$M=N=P$}} & \multicolumn{1}{c}{{Method}} & Relative error& CPU-time (sec) \\ 
				\hline 
				\multirow{4}{*}{$\tt Random$}&\multirow{2}{3em}{$10$}&Algorithm \ref{TG-GMRES} &$3.83\times10^{-2}$&$1.23$\\
				&&Algorithm \ref{TG-GKB}&$1.63\times10^{-10}$& $0.98$\\
				\cline{2-5}
				&\multirow{2}{3em}{$20$}&Algorithm \ref{TG-GMRES}&$3.38\times10^{-2}$&$12.65$\\
				&	&Algorithm \ref{TG-GKB}&$3.12\times10^{-11}$&$6.29$\\
				\hline 
				\multirow{4}{*}{$\tt Uniform$}&\multirow{2}{3em}{$10$}&Algorithm \ref{TG-GMRES} &$4.45\times10^{-2}$&$1.22$\\
				&&Algorithm \ref{TG-GKB}&$3.54\times10^{-6}$& $1.24$\\
				\cline{2-5}
				&\multirow{2}{3em}{$20$}&Algorithm \ref{TG-GMRES}&$1.83\times10^{-1}$&$12.51$\\
				&	&Algorithm \ref{TG-GKB}&$5.05\times10^{-7}$&$10.80$\\
				\hline
				\multirow{4}{*}{$\tt Halton$}&\multirow{2}{3em}{$10$}&Algorithm \ref{TG-GMRES} &$1.96\times10^{-2}$&$1.17$\\
				&&Algorithm \ref{TG-GKB}&$5.73\times10^{-9}$& $1.28$\\
				\cline{2-5}
				&\multirow{2}{3em}{$20$}&Algorithm \ref{TG-GMRES}&$2.77\times10^{-2}$&$12.89$\\
				&	&Algorithm \ref{TG-GKB}&$4.55\times10^{-9}$&$17.09$\\
				\hline  
			\end{tabular}
			\caption{Results for Example 1.}\label{tab1}
		\end{center}
	\end{table}
	\begin{figure}
		\begin{center}
			\includegraphics[width=.33\textwidth,height=.33\textwidth]{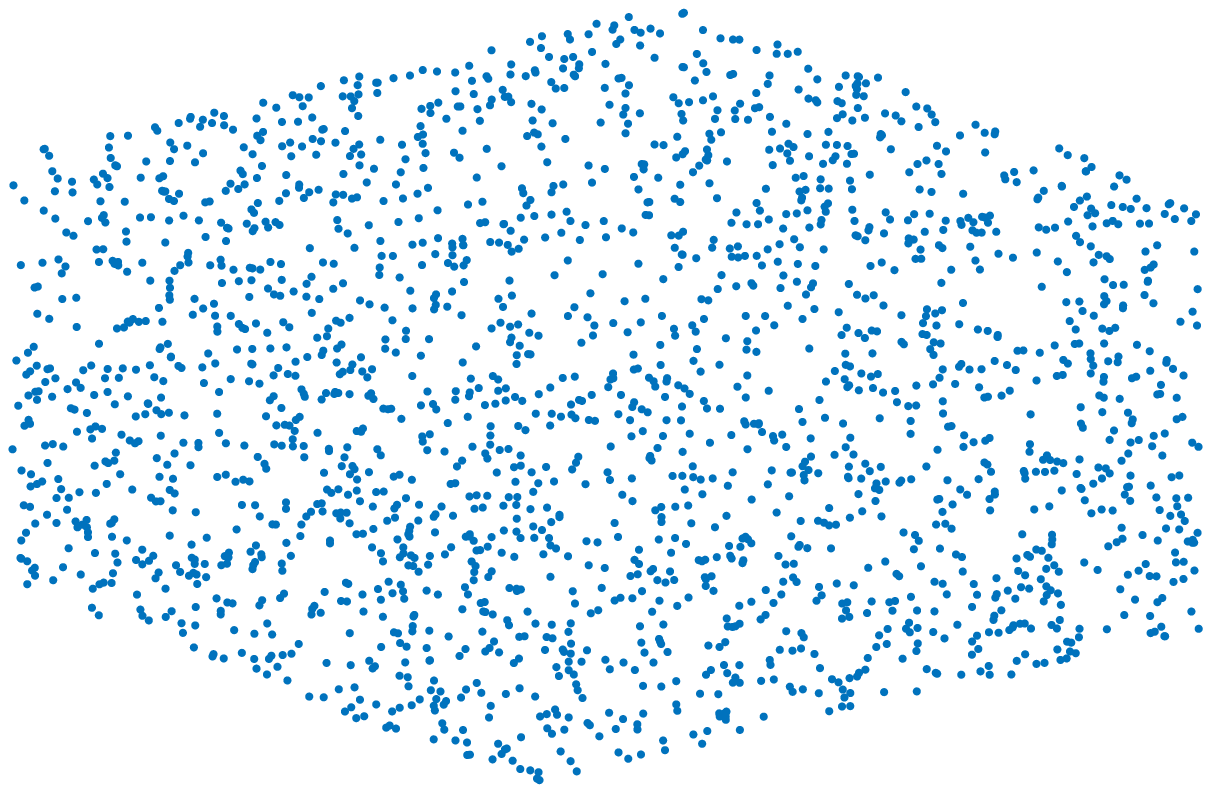}%
			\includegraphics[width=.33\textwidth,height=.33\textwidth]{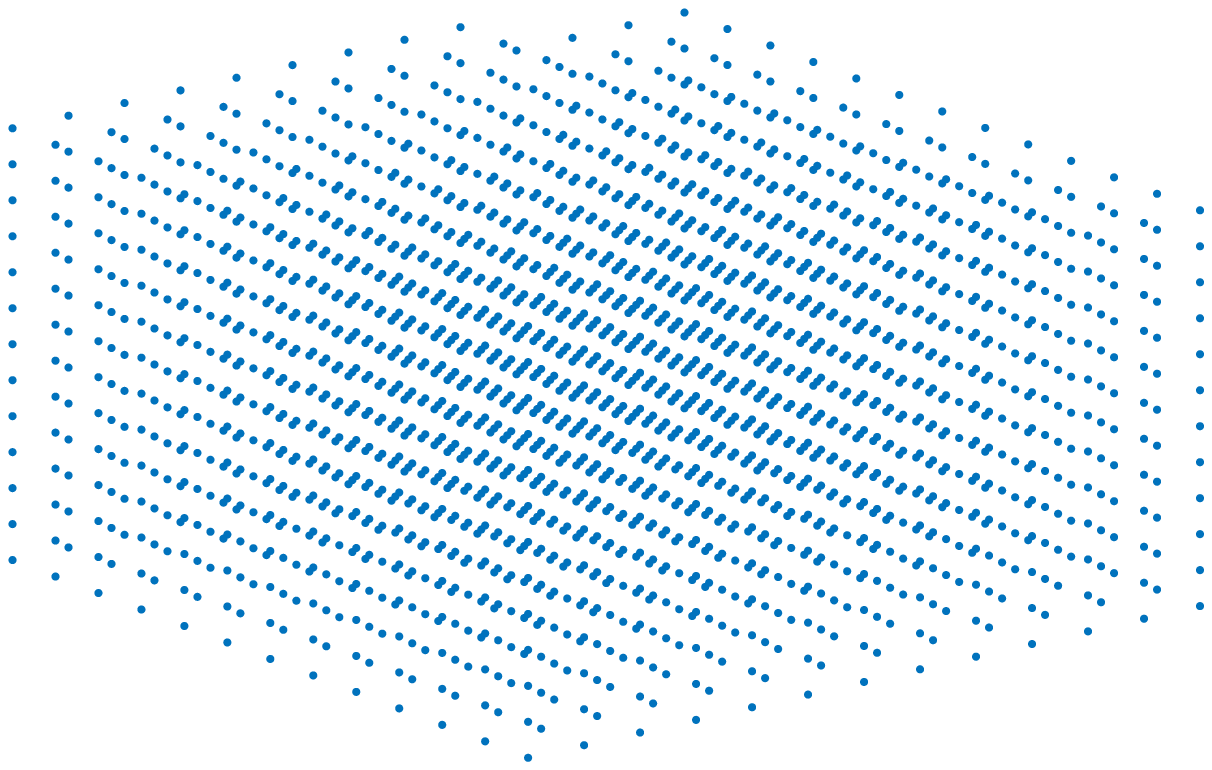}%
			\includegraphics[width=.33\textwidth,height=.33\textwidth]{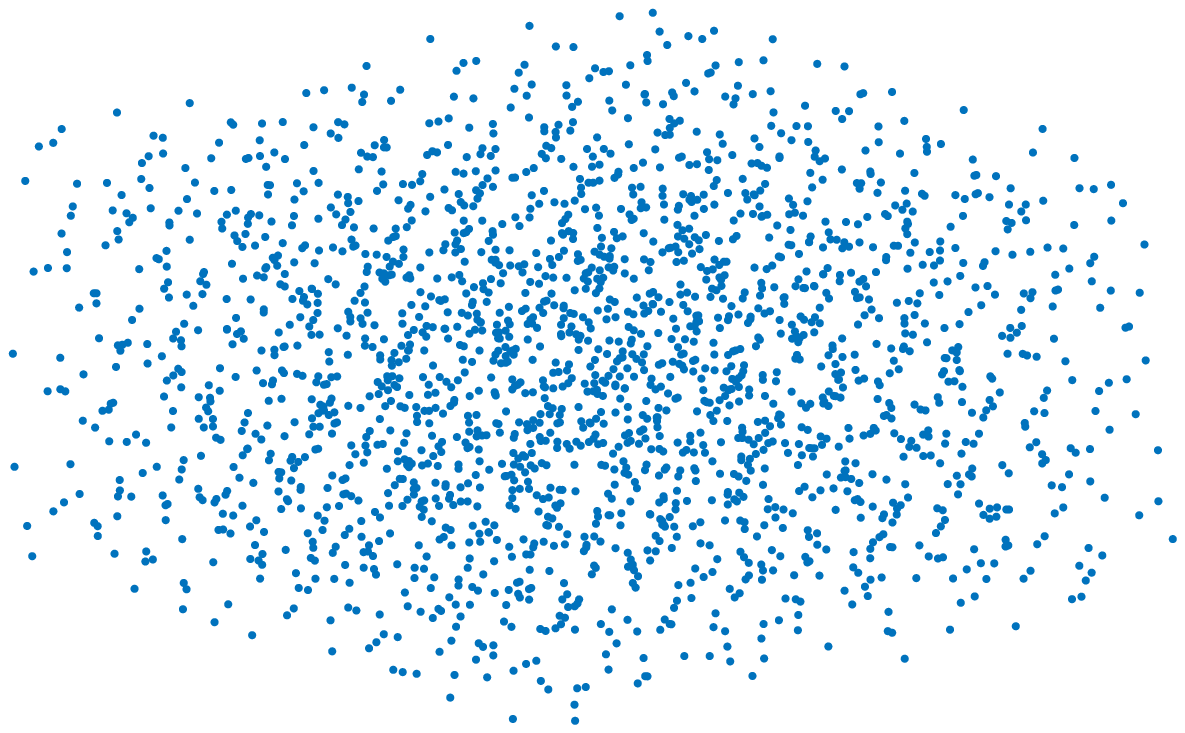}
			\caption{ Sets of collocation points:  random (left),  uniform (middle) and  Halton (right).} \label{fig3}
		\end{center}
	\end{figure} 
	\begin{figure}
		\begin{center}
			\includegraphics[width=.5\textwidth,height=.4\textwidth]{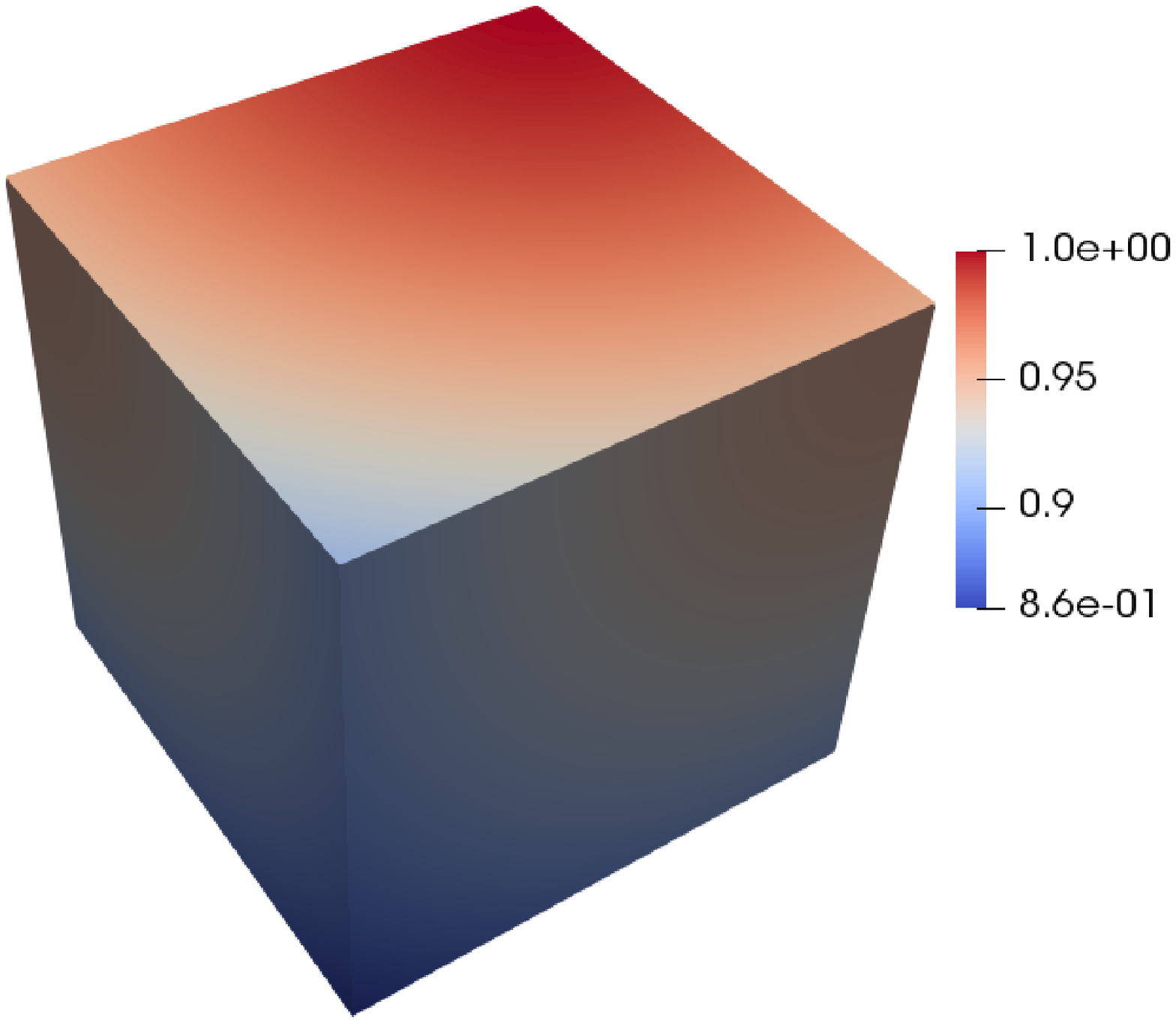}%
			\includegraphics[width=.5\textwidth,height=.4\textwidth]{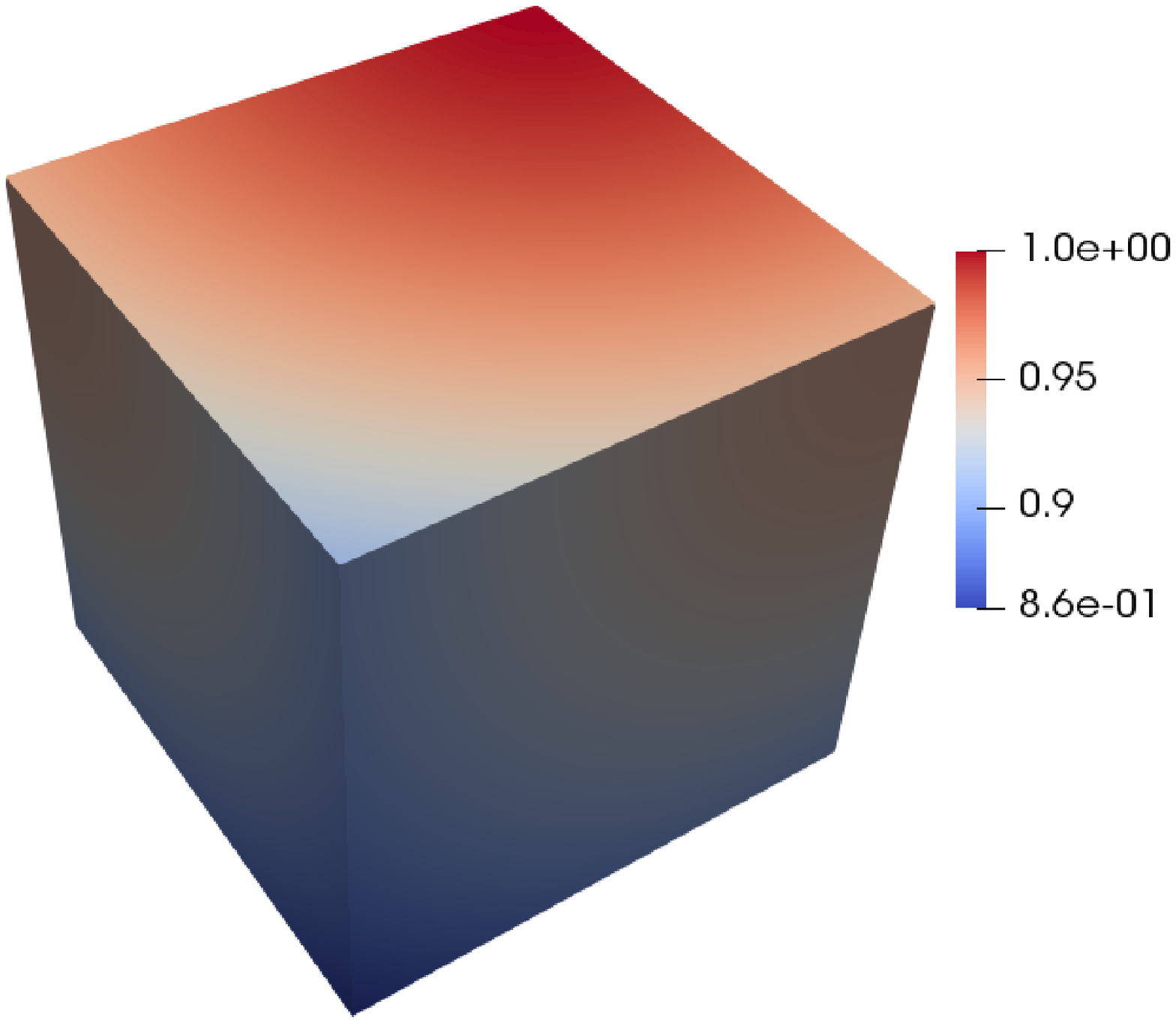}
			\caption{  Exact solution (left) and approximate solution (right).} \label{fig4}
		\end{center}
	\end{figure} 
	\subsection{Example 2} 
	In this example we are interested in the numerical solution of  the three-dimensional Helmholtz equation (\ref{Helm}) in the unit cube $\Omega=[0,1]^3$. This geometry is displayed in Figure \ref{fig1}. As a test the function $f$ is specified so that the exact solution is
	\begin{equation}
	u(x,y,z)=\text{exp}\left(-\frac{x^2+y^2+z^2}{\sigma}\right),
	\end{equation}
	with $\sigma=20$. Table \ref{tab2} displays the performance of Algorithm \ref{TG-GMRES} and Algorithm \ref{TG-GKB}. In Algorithm \ref{TG-GMRES}, we have used as an input for noise level $\nu=10^{-3}$,  $\mathscr{C}$, $\mathscr{X}_0=\mathscr{O}$, $tol=10^{-6}$, $m=10$ and $\text{Iter}_{\text{max}}=10$. The chosen inner and outer iterations  were $m=4$ and  $\text{Iter}_{\text{max}}=4$, respectively. For the ten outer iterations,  minimizing the GCV function  produces  $\mu_{10}=1.15 \times 10^{-5}$. 
	\begin{table}[htbp]
		\begin{center}\footnotesize
			\renewcommand{\arraystretch}{1.3}
			\begin{tabular}{cccccc}\hline
				\multicolumn{1}{c}{{Collocation points}} &\multicolumn{1}{c}{{$M=N=P$}} & \multicolumn{1}{c}{{Method}} & Relative error& CPU-time (sec) \\ 
				\hline 
				\multirow{4}{*}{$\tt Random$}&\multirow{2}{3em}{$10$}&Algorithm \ref{TG-GMRES} &$4.59\times10^{-2}$&$1.15$\\
				&&Algorithm \ref{TG-GKB}&$2.64\times10^{-10}$& $0.73$\\
				\cline{2-5}
				&\multirow{2}{3em}{$20$}&Algorithm \ref{TG-GMRES}&$4.12\times10^{-2}$&$13.22$\\
				&	&Algorithm \ref{TG-GKB}&$5.12\times10^{-11}$&$6.02$\\
				\hline 
				\multirow{4}{*}{$\tt Uniform$}&\multirow{2}{3em}{$10$}&Algorithm \ref{TG-GMRES} &$4.11\times10^{-2}$&$0.22$\\
				&&Algorithm \ref{TG-GKB}&$1.15\times10^{-6}$& $1.64$\\
				\cline{2-5}
				&\multirow{2}{3em}{$20$}&Algorithm \ref{TG-GMRES}&$5.69\times10^{-1}$&$14.65$\\
				&	&Algorithm \ref{TG-GKB}&$2.17\times10^{-6}$&$4.80$\\
				\hline
				\multirow{4}{*}{$\tt Halton$}&\multirow{2}{3em}{$10$}&Algorithm \ref{TG-GMRES} &$1.17\times10^{-1}$&$1.15$\\
				&&Algorithm \ref{TG-GKB}&$4.49\times10^{-5}$& $1.28$\\
				\cline{2-5}
				&\multirow{2}{3em}{$20$}&Algorithm \ref{TG-GMRES}&$1.60\times10^{-1}$&$13.65$\\
				&	&Algorithm \ref{TG-GKB}&$4.50\times10^{-5}$&$18.19$\\
				\hline  
			\end{tabular}
			\caption{Results for Example 2.}\label{tab2}
		\end{center}
	\end{table}
	\subsection{Example 3}
\begin{figure}
	\begin{center}
		\includegraphics[width=.5\textwidth,height=.38\textwidth]{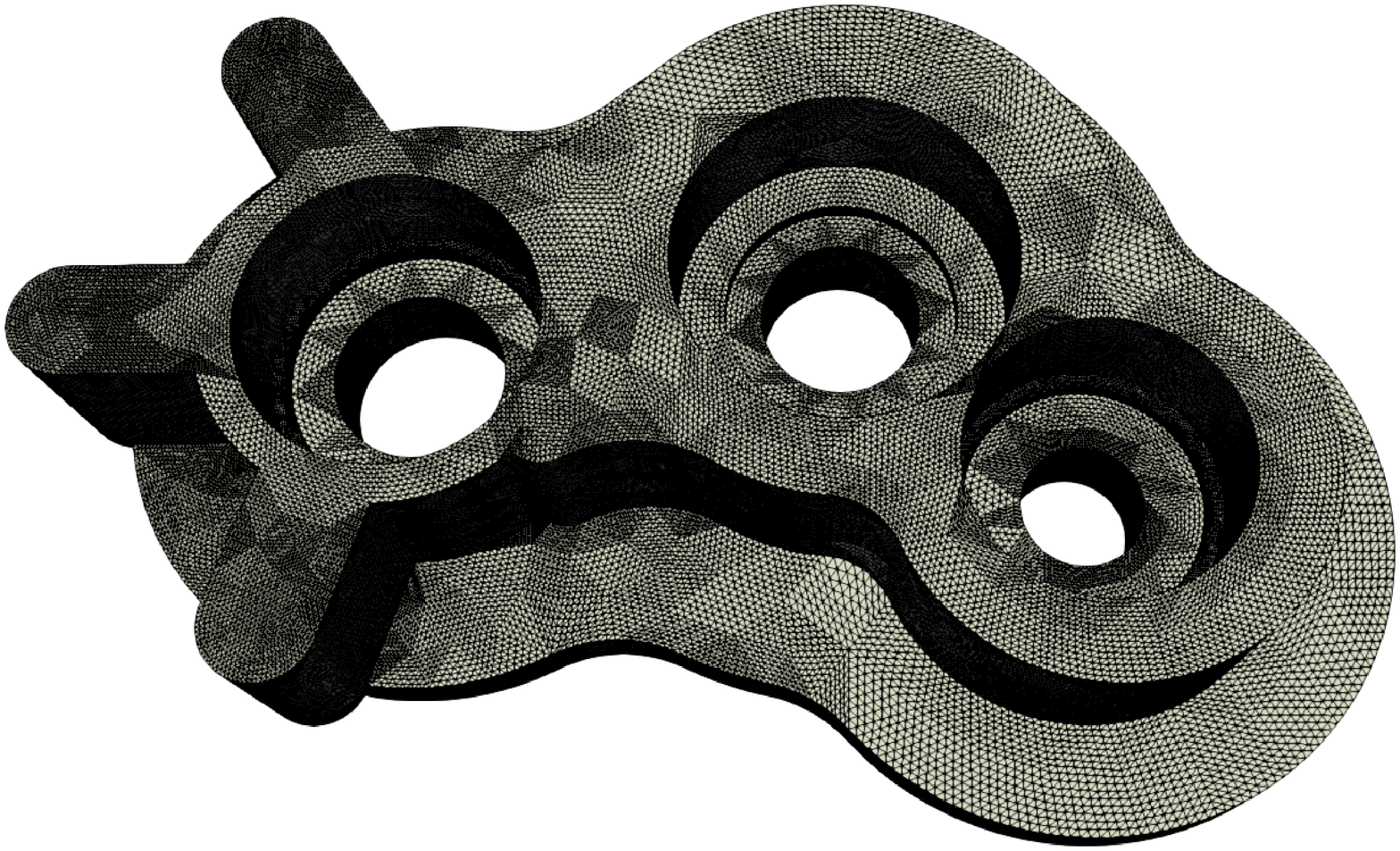}%
		\includegraphics[width=.5\textwidth,height=.4\textwidth]{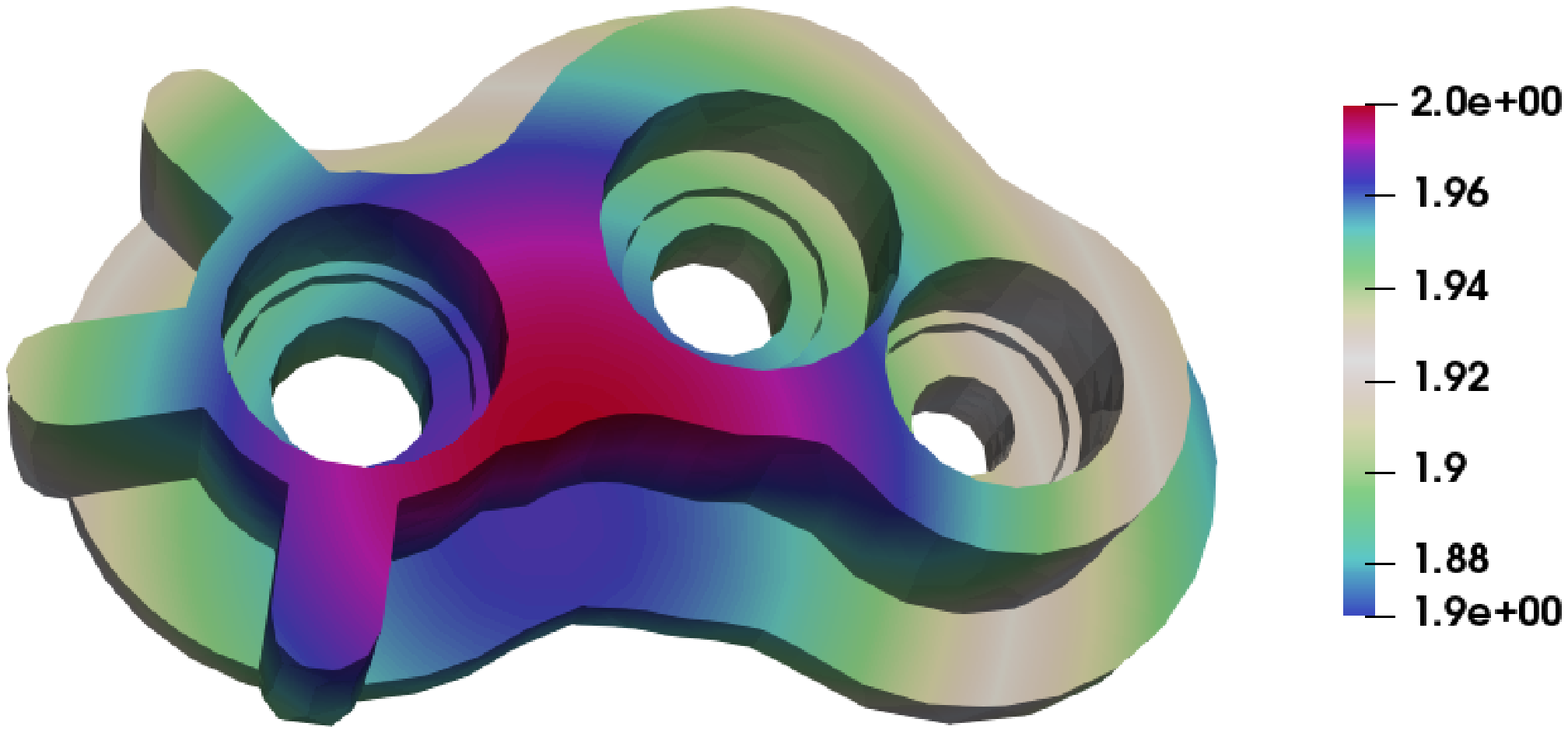}
		\caption{Geometry of the thermal model of a pump casing (left) and exact solution (right).} \label{fig5}
	\end{center}
\end{figure}

\begin{figure}
	\begin{center}
		
		\includegraphics[width=.5\textwidth,height=.4\textwidth]{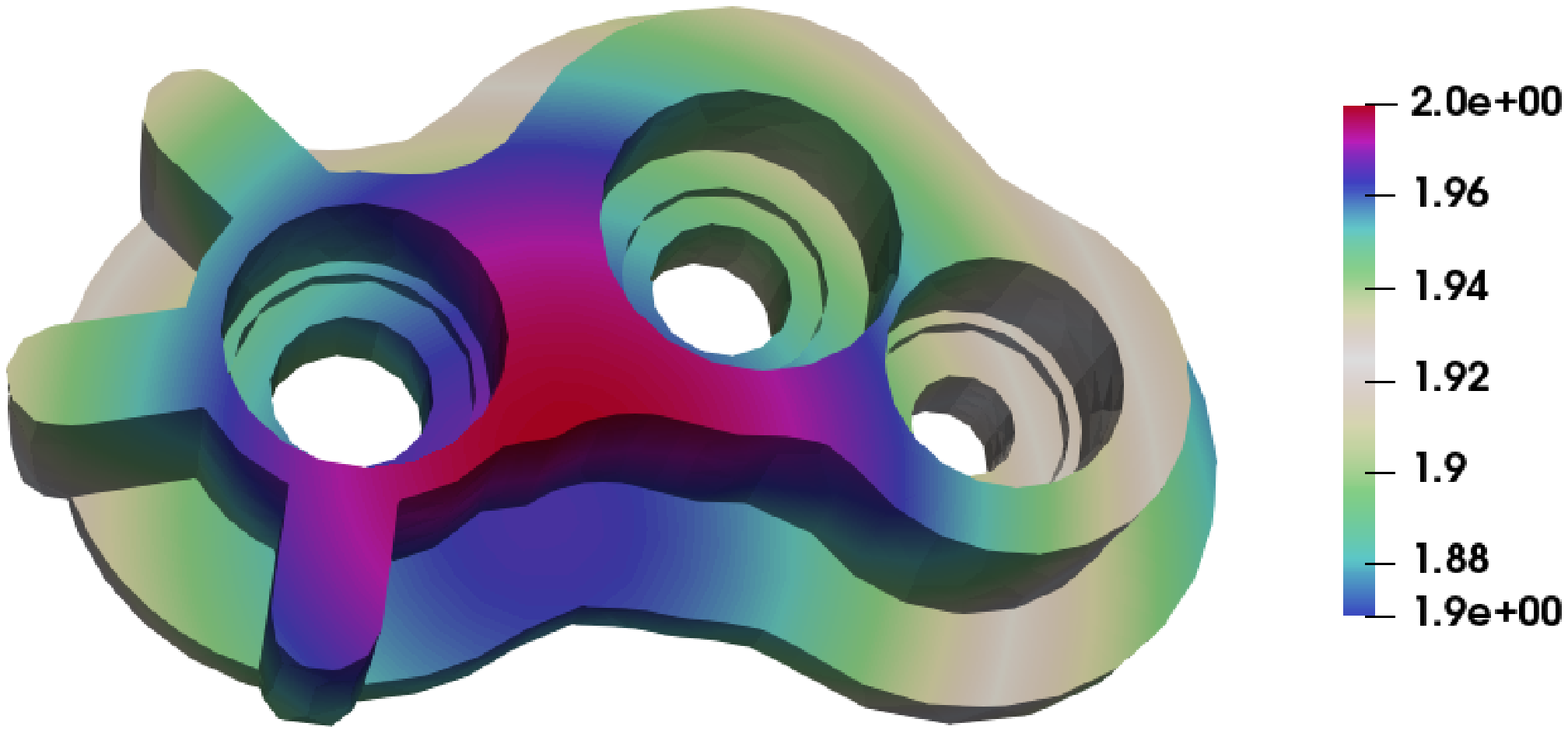}
		\caption{Approximate solution obtained by  Algorithm \ref{TG-GKB}.} \label{fig6}
	\end{center}
\end{figure}

In this example, we illustrate the efficiency of  Algorithm \ref{TG-GKB} applied to solving  ill-posed tensor problem (\ref{tensorprob})  resulting from MQ-RBF discretization of a real-world problem of industrial relevance. We consider the geometry corresponding to a a pump casing model created by using the Gmsh tool \cite{GR}.
Several methods have been proposed in the literature to comprehensively study the
acoustic behaviors of the pump casing \cite{VG, TWGH}. The boundary of the pump model
displayed in Figure \ref{fig5}.  We consider 432124 collocation points. Problems of such large size
add another level of difficulty to our methods, for example, we are unable
to store the underlying tensors $\mathcal{ A}$ and $\mathcal{ H}$  in memory. To overcome this
we resort to the hierarchical MQ-RBF interpolation based compression technique described in Section 4. Specifically. Here,
we consider Dirichlet boundary conditions,
i.e., (\ref{Helmboun}) with $a(x)=1$ and $b(x)=0 ,$. The analytical solution is given by
\begin{equation}
u(x,y,z)=\text{exp}\left(-\frac{(x-1.75)^2+y^2+(z-0.1)^2}{\sigma_1}\right)+\text{exp}\left(-\frac{x^2+y^2+z^2}{\sigma_2}\right),
\end{equation}
with $\sigma_1=\sigma_2=0.5$.  Algorithm \ref{TG-GKB} terminates after 44 steps of Einstein Tensor Global Golub Kahan algorithm to reach $\tau=10^{-12}$, and the computed approximate solution $\mathcal{ X}_{44}$, defined by (\ref{Xkmu}). The relative error corresponding to the computed solution is given by $\frac{\left\|{ \mathcal U}-{\mathcal U}_{\text{computed}}\right\|_{F}}{\|{ \mathcal U}\|_{F}}=1.32\times10^{-9}.$ The computed solution is
shown in Figure \ref{fig6}.
	\section{Conclusion} 
	In this paper we have proposed tensor version of GMRES and Golub–Kahan bidiagonalization
	algorithms using the T-product, with applications to  solving  large-scale linear tensor equations arising in the reconstructions of blurred and noisy multichannel images and videos. The numerical experiments that we have performed  show the
	effectiveness of the proposed schemes to inexpensively computing regularized solutions of high quality.
\bibliographystyle{siam}
\bibliography{Tensor_RBF_3D_2020_v6.bbl}  
\end{document}